\documentclass[a4paper,reqno,11pt]{amsart}
\usepackage{amssymb, amsfonts}
\usepackage{fancyhdr, mathrsfs, graphicx, epsfig, amsmath, amssymb, amsthm}
\usepackage{float}
\usepackage{enumerate}
\usepackage{amscd}
\usepackage{tikz, multicol}
\usetikzlibrary{arrows,decorations.pathmorphing}
\usetikzlibrary{chains}

\theoremstyle{plain}
 \newtheorem{theorem}{Theorem}[section]
 \newtheorem{proposition}[theorem]{Proposition}
 \newtheorem{lemma}[theorem]{Lemma}
 \newtheorem{corollary}[theorem]{Corollary}
 \newtheorem{proposed problem}[theorem]{Proposed problem}
 \newtheorem{open problem}[theorem]{Open problem}
 \newtheorem{proposed project}[theorem]{Proposed project}
\theoremstyle{definition}
 \newtheorem{example}[theorem]{Example}
 \newtheorem{definition}[theorem]{Definition}
\theoremstyle{remark}
 \newtheorem{remark}[theorem]{Remark}
 
 \numberwithin{equation}{section}

\title{Frobenius algebras and root systems: the trigonometric case}
\author{Dali Shen}
\date{\today}
\email{dali\_shen@hotmail.com}

\begin{document}

\setcounter{page}{1}
\pagenumbering{arabic}

\maketitle

\begin{abstract}
We construct Frobenius structures on the $\mathbb{C}^{\times}$-bundle of the complement of a toric
arrangement associated with a root system, by making use of a one-parameter family of torsion free and
flat connections on it. This gives rise to a
trigonometric version of Frobenius algebras in terms of root systems and a new class of Frobenius manifolds.
We also determine their potential functions.
\end{abstract}

\tableofcontents

\section{Introduction}

In this paper, starting from the complement of a toric arrangement associated with a root system,
we construct a Frobenius structure on its $\mathbb{C}^{\times}$-bundle.
This gives rise to a trigonometric version of Frobenius algebras in terms of root systems and
a new class of Frobenius manifolds (in a weak sense, c.f. Definition $\ref{def:Frobenius-manifold}$) as well.

We start from an algebraic torus $H:=\mathrm{Hom}(Q,\mathbb{C}^{\times})$ defined by a root lattice
$Q:=\mathbb{Z}R$ where $R$ is an irreducible reduced root system. We denote the Lie algebra of $H$ by $\mathfrak{h}$
and the Weyl group of this root system $R$ by $W$. Each root $\alpha\in R$ determines a
character $e^{\alpha}:H\rightarrow \mathbb{C}^{\times}$ and gives a corresponding hypertorus
$H_{\alpha}=\{h\in H\mid e^{\alpha}(h)=1\}$. Note that a root $\alpha$ and its negative $-\alpha$ determine
the same $H_{\alpha}$. The collection of these hypertori is called a
toric arrangement associated with a root system, sometimes also referred as a toric mirror arrangement.
We denote by $H^{\circ}$ the complement of this
toric mirror arrangement, i.e.,
\[
H^{\circ}:=H-\cup_{\alpha >0} H_{\alpha}.
\]
Inspired by the work
of Heckman and Opdam \cite{Heckman-Opdam-1, Heckman-Opdam-2}\cite{Opdam-1, Opdam-2} on special hypergeometric
functions associated with root systems, the author constructed a family of torsion free and
flat connections $\tilde{\nabla}^{\kappa}$ on $H^{\circ}\times\mathbb{C}^{\times}$ in \cite{Shen}, depending on the
multiplicity parameter
\[
\kappa=(k_{\alpha})_{\alpha\in R}\in \mathbb{C}^{R},
\]
for which we require it to be a $W$-invariant function so that the ultimately
constructed structure is $W$-invariant.
Since the torsion free and flat connection $\tilde{\nabla}^{\kappa}$ defines an affine structure on
$H^{\circ}\times\mathbb{C}^{\times}$, we would naturally speculate if there exists a Frobenius structure on it.
It turns out to be the case, which is the main theme of this paper.

Taking cue from this torsion free and flat connection $\tilde{\nabla}^{\kappa}$,
we define a product structure for each $\kappa$ on the
tangent bundle of $H^{\circ}\times\mathbb{C}^{\times}$ by
\begin{align*}
  \tilde{X}\cdot_{\kappa}\tilde{Y}
  :=&\frac{1}{2}\sum_{\alpha >0}k_{\alpha}\frac{1+e^{\alpha}}{1-e^{\alpha}}
  \alpha(X) \alpha(Y) \alpha^{\vee}
  -b^{\kappa}(X,Y)-a^{\kappa}(X,Y)t\frac{\partial}{\partial t} \notag \\
&+\lambda_{2}X+\lambda_{1}Y+\lambda_{1}\lambda_{2}t\frac{\partial}{\partial t},
\end{align*}
where $\tilde{X}=X+\lambda_{1} t\frac{\partial}{\partial t}$ is a vector field on $H^{\circ}\times\mathbb{C}^{\times}$
(so is $\tilde{Y}$). The two maps
\[
a^{\kappa}:\mathfrak{h}\times\mathfrak{h}\rightarrow \mathbb{C} \quad \text{and} \quad
b^{\kappa}:\mathfrak{h}\times\mathfrak{h}\rightarrow \mathfrak{h}
\]
are a $W$-invariant symmetric bilinear form and a $W$-equivariant symmetric bilinear map respectively.
More precise definitions for
these notions can be found in Section $\ref{sec:Frobenius-algebra}$.

We can understand the Frobenius structure on the complex manifold by
making use of the so-called structure connection method \cite{Manin},
i.e., a one-parameter family of torsion free and flat connections.
By this we can prove the main theorem in this paper as follows
\begin{theorem}
The product structure $\cdot_{\kappa}$ defined as above on $T(H^{\circ}\times\mathbb{C}^{\times})$
endows each fiber of $T(H^{\circ}\times\mathbb{C}^{\times})$ with a Frobenius algebra structure.
\end{theorem}

We can easily find the identity field $t\frac{\partial}{\partial t}$ for this algebra, then we naturally
have a following corollary.
\begin{corollary}
The manifold $H^{\circ}\times\mathbb{C}^{\times}$ endowed with the structure $(\cdot_{\kappa},a^{\kappa},
t\frac{\partial}{\partial t})$ is a Frobenius manifold.
\end{corollary}

Note that the Euler field is not considered in our definition for Frobenius manifolds (Def. \ref{def:Frobenius-manifold}).

Meanwhile we also have the potential function for this Frobenius structure as follows
\begin{align*}
\Phi=-\frac{t^{3}}{3!}+\frac{t}{2}c^{\kappa}\sum_{\alpha>0}\alpha^{2}
+\sum_{\alpha>0}k_{\alpha}\frac{a^{\kappa}(\alpha^{\vee},\alpha^{\vee})}{\alpha(\alpha^{\vee})}q(\alpha)
-d^{\kappa}\sum_{\alpha>0}\frac{\alpha^{3}}{3!},
\end{align*}
where $t$ stands for the coordinate of the vertical direction, the function $q(\alpha)$ is a series satisfying
\[q'''(\alpha)=\frac{1}{2}\cdot\frac{1+e^{\alpha}}{1-e^{\alpha}},\] the constants $c^{\kappa}$ and $d^{\kappa}$
correspond to the symmetric bilinear form $a^{\kappa}$ and the symmetric cubic form
$a^{\kappa}(b^{\kappa}(\cdot,\cdot),\cdot)$ respectively.

Because of the relation of Frobenius structure with quantum cohomology, it is not surprising that the Frobenius structure
constructed in this paper has a similar form with the work of Bryan and Gholampour \cite{Bryan-Gholampour} on
quantum cohomology of $ADE$ resolutions. By the construction itself, i.e., toric case, one can naturally expect
the potential function would be closely related to the
trigonometric solutions of WDVV equations by Feigin \cite{Feigin}. The configuration in the total space might be
interpreted as an extended $\bigvee$-system due to Stedman and Strachan \cite{Stedman-Strachan}.

The paper is organized as follows. We briefly introduce the definition of Frobenius manifolds in Section
$\ref{intro:Frobenius-manifolds}$, and then construct a Frobenius structure
on our $H^{\circ}\times\mathbb{C}^{\times}$ in Section $\ref{sec:Frobenius-algebra}$, finally in Section
$\ref{sec:toric-Lauricella-manifolds}$ we discuss a class of examples: toric Lauriella manifolds,
which descends to our case when all the weights $\mu_{i}$'s are equal.

\vspace{2mm}
\textbf{Acknowledgements.} I would like to thank Eduard Looijenga for pointing me out this possible direction,
and Di Yang for helpful discussions.
I also thank Shanghai Center for Mathematical Sciences for their generous host during the fall of 2018
where part of this work was done.

\section{Frobenius manifolds and the structure connection}\label{intro:Frobenius-manifolds}

In this section we give a brief introduction to the Frobenius structures on a complex manifold.
For a complete and detailed exposition on this topic,
interested reader can consult the book by Manin \cite{Manin} or Dubrovin \cite{Dubrovin-1996}, as well as
the lecture notes of Looijenga \cite{Looijenga-2009}.

For the moment, for us a $\mathbb{C}$-algebra is simply a $\mathbb{C}$-vector space $A$ endowed with a $\mathbb{C}$-bilinear
map (also referred as the product):
\[
A\times A\rightarrow A;\;(u,v)\mapsto uv
\]
which is associative and a unit element
$e\in A$ such that $e.u=1.u=u$ for all $u\in A$. We often write $1$ for $e$.

\begin{definition}
Let $A$ be a $\mathbb{C}$-algebra which is commutative, associative and finite dimensional as a $\mathbb{C}$-vector space.
A linear function on $A$, $$F:A\rightarrow \mathbb{C}$$ is called a \emph{trace map}
if the map $$a: A\times A \rightarrow \mathbb{C}; \; (u,v) \mapsto
a(u,v):=F(uv)$$ is a nondegenerate bilinear form. The pair $(A,F)$ is called a \emph{Frobenius algebra}.
The bilinear form $a$ sometimes is also called a \emph{pseudometric}.
\end{definition}

\begin{remark}
The fact that the bilinear form $a$ is nondegenerate is equivalent to that the resulting map $u\mapsto F(u.-)$ is a
linear isomorphism of $A$ onto its dual space $A^{*}$ consisting of all the linear forms on $A$. We also need to point
out that the trace map defined here is in general not the one that we usually associate a linear operator (if an
element $u$ of $A$ is regarded as a linear operator $x\in A\mapsto ux\in A$) with its trace.
\end{remark}

\begin{lemma}\label{lem:associativity-law}
The bilinear form $a$ satisfies the associative law $a(uv,w)=a(u,vw)$. And conversely, any nondegenerate bilinear
symmetric map $a:A\times A\rightarrow \mathbb{C}$ with the associative law determines a trace map on $A$.
\end{lemma}

\begin{proof}
$a(uv,w)=F((uv)w)=F(u(vw))=a(u,vw)$ since $A$ is an associative $\mathbb{C}$-algebra, then the first statement follows.

Conversely, we can define a linear function $I$ by $I(u):=a(u,e)$. Then we can define a new map $a':A\times A\rightarrow\mathbb{C}$
as follows: $a'(u,v):=I(uv)$, but we have $I(uv)=a(uv,e)=a(u,v)$ by the associativity of $a$. This shows the newly
defined map $a'$ is the same as $a$, which is also a nondegenerate bilinear symmetric form. The second statement follows.
\end{proof}

This associativity law of the bilinear form $a$ is also called a \emph{Frobenius condition}.

Here are some simple examples of Frobenius algebra.

\begin{example}
(i) For the field $\mathbb{C}$ which could be viewed as a $\mathbb{C}$-algebra, we can define a linear form
by a nonzero scalar multiplication $F:\mathbb{C}\rightarrow \mathbb{C};\; u\mapsto \nu u$ for $\nu\neq 0$.
This is a trace map on $\mathbb{C}$.

(ii) Let $A=\mathbb{C}[t]/(t^{n})$ with $n\in\mathbb{Z}_{+}$. A linear form $F:A\rightarrow \mathbb{C}$ is a trace map
if and only if $F(t^{n-1})\neq 0$.
\end{example}

Now let us see what kind of role the associativity condition plays here? If we are given a Frobenius algebra, we can
also consider the trilinear map $T: A\times A\times A\rightarrow \mathbb{C}$ defined by $T(u,v,w):=F(uvw)$. But conversely,
if we are only given a vector space $A$, a trilinear map $T:A\times A\times A\rightarrow \mathbb{C}$, and an element
$e\in A$, does $T$ define a Frobenius algebra structure on $A$? The answer is obviously no. We must impose some additional
conditions so that $T$ can be used to define a Frobenius algebra on $A$. First $T$ must be required to be symmetric.
And we also want that the bilinear form $(u,v)\in A \times A\mapsto T(u,v,e)\in\mathbb{C}$ is nondegenerate. We thus
have defined a bilinear map (product) $A \times A\rightarrow A; \; (u,v)\mapsto uv$ characterized by that $T(uv,x,e)=T(u,v,x)$
for all $x\in A$. Since $T$ is symmetric, the product is commutative. And $e$ becomes the identity element of $A$ for
$ue$ is characterized by $T(ue,e,x)=T(u,e,x)$ for all $x\in A$ which implies $ue=u$.

Besides these two conditions, the associativity
does not hold a priori and thus has to be endowed. This means we want that $T(uv,w,x)=T(u,vw,x)$ for all $u,v,w,x\in A$.
In fact, we can write out this condition in terms of a basis of $A$. If $\{u_{1},\cdots,u_{n}\}$ is a basis of $A$, define
$T_{ijk}:=T(u_{i},u_{j},u_{k})$, then $(a_{jk}:=T_{1jk})_{jk}$ is a nondegenerate matrix. Let $(a^{jk})_{jk}$ denote
its inverse matrix, then by the above recipe: $T(u_{i},u_{j},u_{k})=T(u_{i}u_{j},u_{k},e)=T(b^{l}u_{l},u_{k},e)$, we have
\[ u_{i}u_{j}=T_{ijk}a^{kl}u_{l}  \]
where the Einstein summation convention is used. The above associativity condition also means we want that
$T(u_{i}u_{j},u_{k},u_{l})=T(u_{i},u_{j}u_{k},u_{l})$ in terms of the basis. So we have
$T(T_{ijp}a^{pq}u_{q},u_{k},u_{l})=T(u_{i},T_{jkp}a^{pq}u_{q},u_{l})$ which is equivalent to
\begin{align*}\label{Associativity}\tag{Ass.}
T_{ijp}a^{pq}T_{qkl}=T_{jkp}a^{pq}T_{iql}.
\end{align*}
This is a system of equations which must be satisfied in order that the product being associative.

\vspace{1mm}
Now let be given a complex manifold $M$ whose holomorphic tangent bundle is denoted by $TM$. We are also given on $TM$
a nondegenerate symmetric bilinear form $a$ and a symmetric trilinear form $T$, both depending holomorphically on the
base point. The product of this bundle can be characterized by the property that $a(XY,Z)=T(X,Y,Z)$, denoted by
\[
\cdot:TM\times TM\rightarrow TM; \; X\cdot Y\mapsto XY.
\]
It is clear that this product is commutative by the symmetry
of $T$. We use $\nabla$ to denote the complex counterpart of the $\emph{Levi-Civita}$ connection on the holomorphic
tangent bundle $TM$ which is characterized by the following $2$ properties:
\leftmargini=7mm
\begin{description}
\item [compatibility] $Z(a(X,Y))=a(\nabla_{Z}X,Y)+a(X,\nabla_{Z}Y)$,
\item [torsion freeness] $\nabla_{X}Y-\nabla_{Y}X=[X,Y]$.
\end{description}
Its curvature form is given by
\[ \mathrm{R}(\nabla)(X,Y)Z:=\nabla_{X}\nabla_{Y}Z-\nabla_{Y}\nabla_{X}Z-\nabla_{[X,Y]}Z. \]

We can then define a one-parameter family of connections $\nabla(\mu)$ on this bundle by
\[ \nabla(\mu)_{X}Y:=\nabla_{X}Y+\mu X\cdot Y, \quad \mu\in \mathbb{C}. \]
The connection $\nabla(\mu)$ is called the \emph{structure connection} of $(M,a,T)$.

By the commutativity of the product we immediately have
\[ \nabla(\mu)_{X}Y-\nabla(\mu)_{Y}X-[X,Y]=\nabla_{X}Y-\nabla_{Y}X-[X,Y]=0, \]
which shows that $\nabla(\mu)$ is torsion free as well.
If for a local vector field $X$ on $M$, $\iota_{X}$ denotes the multiplication
operator on vector fields:
\[\iota_{X}(Y):=X\cdot Y, \]
then we can define a new tensor
\[ \mathrm{R}'(\nabla)(X,Y):=[\nabla_{X},\iota_{Y}]-[\nabla_{Y},\iota_{X}]-\iota_{[X,Y]} \]
which is a holomorphic $2$-form taking values in the symmetric endomorphism of $TM$. It's clear that
$a(\mathrm{R}'(\nabla)(X,Y)Z,W)$ is antisymmetric in $(X,Y)$ and symmetric in $(Z,W)$.

\begin{proposition}\label{Frobenius-conditions}
The following statements are equivalent:
\begin{enumerate}[(i)]
\item $\nabla$ is flat, the product is associative and if $X,Y,Z$ are (local) flat vector fields on a domain $U\subset M$,
then the trilinear form $T(X,Y,Z)$ locally is given by
$T(X,Y,Z)=\nabla_{X}\nabla_{Y}\nabla_{Z}\Phi$ where $\Phi:U\rightarrow \mathbb{C}$ is a holomorphic
function on $U$.
\item $\nabla$ is flat, the product is associative and $\mathrm{R}'\equiv 0$.
\item The connection $\nabla(\mu)$ is flat for any $\mu\in\mathbb{C}$.
\end{enumerate}
\end{proposition}

\begin{proof}
First prove $(ii)\Leftrightarrow (iii)$. We have
\begin{align*}
\nabla(\mu)_{X}\nabla(\mu)_{Y}&=(\nabla_{X}+\mu \iota_{X})(\nabla_{Y}+\mu \iota_{Y})\\
   &=\nabla_{X}\nabla_{Y}+\mu(\iota_{X}\nabla_{Y}+\nabla_{X}\iota_{Y})+\mu^{2}\iota_{X}\iota_{Y}.
\end{align*}
Similarly,
\begin{align*}
\nabla(\mu)_{Y}\nabla(\mu)_{X}&
  =\nabla_{Y}\nabla_{X}+\mu(\iota_{Y}\nabla_{X}+\nabla_{Y}\iota_{X})+\mu^{2}\iota_{Y}\iota_{X},\\
\nabla(\mu)_{[X,Y]}&=\nabla_{[X,Y]}+\mu\iota_{[X,Y]}.
\end{align*}
Then we have
\begin{align*}
\mathrm{R}(\nabla(\mu))(X,Y)&=\nabla(\mu)_{X}\nabla(\mu)_{Y}-\nabla(\mu)_{Y}\nabla(\mu)_{X}-\nabla(\mu)_{[X,Y]}\\
&=\mathrm{R}(\nabla)(X,Y)+\mu \mathrm{R}'(\nabla)(X,Y)+\mu^{2}(\iota_{X}\iota_{Y}-\iota_{Y}\iota_{X}).
\end{align*}
So if $\nabla$ is flat, i.e., $\mathrm{R}(\nabla)=0$. We then see that $\nabla(\mu)$ is flat for all $\mu$ if and
only if $\mathrm{R}'(\nabla)=0$ and $\iota_{X}\iota_{Y}=\iota_{Y}\iota_{X}$ for all $X,Y$. While the condition that
$\iota_{X}\iota_{Y}=\iota_{Y}\iota_{X}$ for all $X,Y$ is equivalent to that $X\cdot(Y\cdot Z)=Y\cdot(X\cdot Z)$ for all
$X,Y,Z$. But by the commutativity of the product the left hand side $X\cdot(Y\cdot Z)=X\cdot(Z\cdot Y)$,
and the right hand side $Y\cdot(X\cdot Z)=(X\cdot Z)\cdot Y$.
This is just the associativity property. So $(ii)\Leftrightarrow (iii)$ follows.

Now let us prove $(i)\Leftrightarrow (ii)$. Since $a$ is flat we can pass all the things to a flat chart $(U,\varphi)$
such that $D=\varphi (U)\subset \mathbb{C}^{n}$ is an open polydisk. Under this setting, $a$ has constant coefficients,
$\nabla$ becomes the usual derivation and the flat vector fields are just the constant ones. Suppose we are given holomorphic
functions $f_{ijk}:D\rightarrow \mathbb{C}$ for $1\leq i,j,k\leq n$. It is well-known that these can arise as the
third order partial derivatives of a holomorphic function $\Phi$ if and only if $\partial_{l}f_{ijk}$ is symmetric in
all its indices. In other words, if $f$ is a trilinear form on the tangent bundle of $D$, then there exists a
holomorphic function $\Phi$ such that $f(X,Y,Z)=\nabla_{X}\nabla_{Y}\nabla_{Z}\Phi$ for all triples of flat vector fields
$(X,Y,Z)$ if and only if $X(f(Y,Z,W))$ is symmetric in its all arguments for all quadruples of flat vector fields
$(X,Y,Z,W)$. Since we have $f(X,Y,Z)=a(XY,Z)$ and we already know that $f$ is symmetric in its three arguments.
We have
\begin{align*}
Xa(Y\cdot Z,W)&=a(\nabla_{X}(Y\cdot Z),W)+a(Y\cdot Z,\nabla_{X}W)\\
  &=a(\nabla_{X}(Y\cdot Z),W).
\end{align*}
But since $X,Y,Z$ are all flat, we also have
\begin{align*}
\mathrm{R}'(\nabla)(X,Y)Z&=\nabla_{X}(Y\cdot Z)-\nabla_{Y}(X\cdot Z),
\end{align*}
and then it is clear that $Xa(Y\cdot Z,W)$ is symmetric in $X$ and $Y$ if and only if $\mathrm{R}'=0$.
\end{proof}

\begin{remark}
If we denote the coefficient of $\mu^{2}$ in $\mathrm{R}(\nabla(\mu))$, i.e., the tensor $\iota_{X}\iota_{Y}-\iota_{Y}\iota_{X}$,
by $\mathrm{R}''(X,Y)$, then from the above proposition we can see that
(i) the condition $\mathrm{R}'=0$ is a potential condition, and
(ii) the condition $\mathrm{R}''=0$ is an associativity condition.
\end{remark}

\begin{remark}
The function $\Phi$ that appears in Statement (i) of Proposition $\ref{Frobenius-conditions}$ is called a
(local) \emph{potential function}. Since here only its third order derivatives matter, it is
(in terms of flat coordinates $(z^{1},\cdots,z^{n})$) unique up to a polynomial of degree two.
In particular, a potential function needs not be defined on all of $M$. The associativity equation
(\ref{Associativity}) now is read as a highly nontrivial system of partial differential equations:
if $(z^{1},\cdots,z^{n})$ is a system of flat coordinates and $\partial_{\nu}:=\frac{\partial}{\partial z^{\nu}}$,
then we require that for all $i,j,k,l$,
\[
(\partial_{i}\partial_{j}\partial_{p}\Phi)a^{pq}(\partial_{q}\partial_{k}\partial_{l}\Phi)=
(\partial_{j}\partial_{k}\partial_{p}\Phi)a^{pq}(\partial_{i}\partial_{q}\partial_{l}\Phi).  \tag{WDVV}
\]
These are known as the \emph{Witten-Dijkgraaf-Verlinde-Verlinde} equations.
\end{remark}

Then we are properly prepared to introduce the main notion of this section.

\begin{definition}\label{def:Frobenius-manifold}
A complex manifold $M$ is called a $\emph{Frobenius manifold}$
if its holomorphic tangent bundle is fiberwisely endowed
with the structure of a Frobenius algebra $(\cdot,F,e)$ satisfying
\begin{enumerate}[(i)]
\item the equivalent conditions of Proposition $\ref{Frobenius-conditions}$ are fulfilled for the associated symmetric
bilinear and trilinear forms $a$ and $T$, \\
\noindent and
\item the identity field $e$ on $M$ is flat for the Levi-Civita connection of $a$.
\end{enumerate}
\end{definition}

\begin{remark}
Note that Dubrovin \cite{Dubrovin-1996} requires an Euler vector field for the definition of a Frobenius manifold.
And Manin in his book \cite{Manin} starts from a $\mathbb{Z}/2\mathbb{Z}$-graded structure sheaf on a manifold
(which he called a supermanifold) for his definition
of a Frobenius manifold. But in this paper we do not introduce these notions because we want to focus on the aforementioned
more central conditions for our construction of Frobenius algebras associated with root systems.
\end{remark}

Here are some examples of Frobenius manifolds.

\begin{example}
(i) The trivial example is $M=\mathbb{C}^{n}$ whose coordinates are $(z^{1},\cdots,z^{n})$, $a=\sum_{i}(dz^{i})^{2}$
and product $\partial_{i}\cdot\partial_{i}=\partial_{i}$. A potential function is a cubic form $\Phi(z)=\frac{1}{6}
\sum_{i}(z^{i})^{3}$ and the family of connections is given by $\nabla(\mu)_{\partial_{i}}\partial_{j}=
\mu\delta_{j}^{i}\partial_{i}$.

(ii) (Two-dimensional case) In this case the product on a vector space $A$ of dimension two with nonzero
unit $e$ is automatically associative. We then have
$A$ is isomorphic to the semisimple $\mathbb{C}\oplus\mathbb{C}$ or to the nonsemisimple $\mathbb{C}[y]/(y^{2})$.
It remains to find the potential functions. Let
$e$ be the unit vector field and $F$ the trace differential. Since $e$ is flat, $a(e,e)=F(e\cdot e)$ is constant,
say equal to $c\in\mathbb{C}$. There are two cases depending on whether $c$ is 0 or not.

We first do the case
$c=0$. Then we can find flat coordinates $(z,w)$ such that $e=\partial_{z}$ and $a=dz\otimes dw+dw\otimes dz$.
Since we have $a(\partial_{z}\cdot\partial_{z},\partial_{z})=a(\partial_{z},\partial_{z})=0$
and $a(\partial_{z}\cdot\partial_{z},\partial_{w})=a(\partial_{z},\partial_{w})=1$, it follows that
$\Phi_{zzz}=0$ and $\Phi_{zzw}=1$. But since $\partial_{z}\cdot\partial_{w}=\partial_{w}$, we must also
have $\Phi_{zwz}=1$ and $\Phi_{zww}=0$. It follows that $\Phi(z,w)=\frac{1}{2}z^{2}w+f(w)$ up to
quadratic terms, where $f$ is holomorphic.

If $c\neq 0$, then we can find flat coordinates $(z,w)$ such that $e=\partial_{z}$ and $a=cdz\otimes dz+cdw\otimes dw$.
Then we want that $\Phi_{zzz}=c$, $\Phi_{zzw}=0$, $\Phi_{zwz}=0$, $\Phi_{zww}=c$. It follows that
$\Phi(z,w)=\frac{1}{6}cz^{3}+\frac{1}{2}czw^{2}+f(w)$ up to quadratic terms, where $f$ is holomorphic.

Conversely, in both cases, with these choices of $e$ and $a$, any $\Phi$ of the form defines a Frobenius manifold.
\end{example}

\begin{remark}
The most important class of examples is furnished by quantum cohomology which in fact motivated the
definition in the first place. And another beautiful class of examples is
furnished by the space of polynomials which is due to Saito \cite{Saito} and Dubrovin \cite{Dubrovin-1992}.
But we will not elaborate these two
important classes of Frobenius manifolds over here. Interested readers can consult Manin \cite{Manin} for
detailed explanation.
\end{remark}

\section{Frobenius algebras from toric Dunkl connection}\label{sec:Frobenius-algebra}

In the previous work of the author \cite{Shen}, starting from a toric mirror arrangement complement,
we have constructed affine structures on its $\mathbb{C}^{\times}$-bundle by showing that there exists
a family of torsion free and flat connections
on this total space.
By regarding them as a structure connection we can thus define a product structure
on the tangent bundle of this manifold. This gives rise to a new class of Frobenius manifolds associated with
root systems.

\subsection{Root systems}

Let $\mathfrak{a}$ be a real vector space of dimension $n$, which is further made to be a Euclidean vector space by
endowing it with an inner product $(\cdot,\cdot)$. Denote its dual vector space by $\mathfrak{a}^{*}$. We can identify
$\mathfrak{a}$ with $\mathfrak{a}^{*}$ by the inner product, so that the dual space $\mathfrak{a}^{*}$ is also
endowed with an inner product, denoted by $(\cdot,\cdot)$ as well by abuse of notation.

For a nonzero vector $\alpha\in \mathfrak{a}^{*}$, there corresponds an orthogonal
reflection $s_{\alpha}$ with the hyperplane perpendicular to $\alpha$
being the mirror. This reflection could be written as
\[ s_{\alpha}(\beta)=\beta-\frac{2(\beta,\alpha)}{(\alpha,\alpha)}\alpha \]
for any $\beta\in \mathfrak{a}^{*}$. We can easily check that
\[ s_{\alpha}(\alpha)=-\alpha \; \text{and} \; s_{\alpha}(\beta)=\beta
\; \text{for} \; (\beta,\alpha)=0. \]
Then $s_{\alpha}^{2}=1$ follows from the above formula directly.
We recall the definition of a root system first.

A finite subset $R$ of $\mathfrak{a}^{*}$ is called a $\emph{root system}$ if it does not contain $0$
and spans $\mathfrak{a}^{*}$ such that any $s_{\alpha}$ leaves $R$ invariant and
$s_{\alpha}(\beta)\in \beta+\mathbb{Z}\alpha$ for any $\alpha,\beta\in R$. Any vector
belonging to $R$ is called a root. The dimension of $\mathfrak{a}^{*}$ is called the
$\emph{rank}$ of the system. The group $W(R)$ generated by the $s_{\alpha}$ is
called the $\emph{Weyl group}$ of $R$.
This root system $R$ is said to be $\emph{reduced}$ if $R\cap\mathbb{R}\alpha
=\{\alpha,-\alpha\}$ for any $\alpha\in R$, and said to be $\emph{irreducible}$ if
nonempty $R$ can not be decomposed as a direct sum of two nonempty root systems.

For each $\alpha\in R$ there exists a coroot $\alpha^{\vee}\in\mathfrak{a}$ such that
$(\alpha,\alpha^{\vee})=2$ and $(\beta,\alpha^{\vee})\in\mathbb{Z}$ for all $\alpha,\beta\in R$, and for any $\alpha\in R$
the reflection $s_{\alpha}(\gamma)=\gamma-(\gamma,\alpha^{\vee})\alpha$ leaves $R$ invariant. The set
$R^{\vee}=\{\alpha^{\vee}|\alpha\in R\}$ is again a root system in $\mathfrak{a}$, called the coroot system relative to
$R$.

\vspace{1mm}
Suppose now we are given a reduced irreducible root system $R \subset \mathfrak{a}^{*}$.
The integral span $Q=\mathbb{Z}R$ of the root system $R$ in $\mathfrak{a}^{*}$ is called the root lattice, its dual
$P^{\vee}=\mathrm{Hom}(Q,\mathbb{Z})$ in $\mathfrak{a}$ is called the coweight lattice of $R^{\vee}$.
Hence we have an algebraic torus defined as follows
\[
H=\mathrm{Hom}(Q,\mathbb{C}^{\times})
\]
with character lattice being $Q$, sometimes also called an adjoint torus.

We denote by $\mathfrak{h}$ the Lie algebra of $H$, which is equal to $\mathbb{C}\otimes P^{\vee}$.
First let us consider a $W$-invariant symmetric bilinear form $a$ and a $W$-equivariant symmetric bilinear map $b$
on $\mathfrak{h}$ respectively as follows
\[
a: \mathfrak{h}\times\mathfrak{h}\rightarrow\mathbb{C},\quad b:\mathfrak{h}\times\mathfrak{h}\rightarrow\mathfrak{h}.
\]
We have the following characterization for $a$ and $b$.

\begin{lemma}\label{bmap}
Let $a$ and $b$ be given as above.
If $R$ is irreducible then
\begin{enumerate}
\item The $W$-invariant symmetric bilinear form $a$ is just a multiple of the given inner product.

\item The $W$-equivariant symmetric bilinear map
$b$ vanishes unless $R$ is of type $A_{n}$ for
$n\geq 2$ in which case there exists a $k'\in\mathbb{C}$ such that
\[ b(u,v)=\frac{1}{2}k'\sum_{\alpha >0}\alpha(u)\alpha(v)\alpha'  \quad
\text{for any} \; u,v\in \mathfrak{h} \]
with $\alpha'=\varepsilon_{i}+\varepsilon_{j}-\frac{2}{n+1}\sum_{l}\varepsilon_{l}$
if we take the construction of $\alpha$ from Bourbaki \cite{Bourbaki}:
$\alpha=z_{i}-z_{j}$ for $1\leq i<j \leq n+1$, where $\{z_{i}\}$ is the dual basis of $\{\varepsilon_{i}\}$
in $(\mathbb{R}^{n+1})^{*}$.
\end{enumerate}
\end{lemma}

\begin{proof}
(1) The given inner product $(\cdot,\cdot)$ on $\mathfrak{a}$ can be extended $\mathbb{C}$-linearly
to a nondegenerate symmetric bilinear form
on $\mathfrak{h}$, which is invariant under the action of $W$, still denoted by $(\cdot,\cdot)$.
Since $R$ is irreducible, then the $W$-invariant
symmetric bilinear form $a$ is just a multiple of this given inner product $(\cdot,\cdot)$ by Schur's lemma.

(2) See Lemma 2.5 of \cite{Shen}.
\end{proof}

\begin{remark}\label{rem:k'-for-An}
In fact, for type $A_{n}$, besides the construction in the proof of \cite{Shen},
there is also another way to obtain a generator for $b$: by taking $v\in\mathfrak{a}
\mapsto\partial_{v}\bar{\sigma}_{3}$, where
$\bar{\sigma}_{3}:=\sigma_{3}|\mathfrak{a}^{*}$ is an element of
$(\mathrm{Sym}^{3}\mathfrak{a}^{*})^{W}$. This viewpoint will become more clear
when we discuss the toric Lauricella case in Section $\ref{sec:toric-Lauricella-manifolds}$.
\end{remark}

\subsection{Frobenius structures}

Each root $\alpha$ of $R$ is primitive in $Q$ and determines a character $e^{\alpha}: H\rightarrow\mathbb{C}^{\times}$.
The kernel of the character $H_{\alpha}=\{h\in H\mid e^{\alpha}(h)=1\}$ defines a hypertorus, called the mirror
determined by $\alpha$. Its Lie algebra $\mathfrak{h}_{\alpha}$ is the zero set of $\alpha$ in $\mathfrak{h}$.
The root system is closed under inversion and note that the negative $-\alpha$ determines the same hypertorus as $\alpha$.
The finite collection of these hypertori $H_{\alpha}$'s
is called a \emph{toric mirror arrangement} associated to
the root system $R$, sometimes in this paper also referred as a \emph{toric arrangement} for short if it leads no confusion.
We write $H^{\circ}$ for the complement of the toric mirror arrangement as follows:
\[
H^{\circ}:=H-\cup_{\alpha>0}H_{\alpha}.
\]

For $u\in\mathfrak{h}$ we denote by $\partial_{u}$ the associated
translation invariant vector field on $H$. Likewise, for $\phi\in\mathfrak{h}^{*}$ we denote by $d\phi$ the
associated translation invariant differential on $H$. In case $\phi\in Q$, it determines a character of $H$,
$e^{\phi}:H\rightarrow\mathbb{C}^{\times}$, then we have $d\phi=(e^{\phi})^{*}(\frac{dt}{t})$ with $t$ the coordinate on
$\mathbb{C}^{\times}$. We denote by $\nabla^{0}$ the flat translation invariant connection on $H$, so that
$\nabla^{0}_{\partial_{v}}=\partial_{v}$. So is $\tilde{\nabla}^{0}$ on $H\times\mathbb{C}^{\times}$.

Let $\kappa$ be a $W$-invariant function
\[
\kappa=(k_{\alpha})_{\alpha\in R}\in\mathbb{C}^{R},
\]
meaning that $k_{w\alpha}=k_{\alpha}$ for any $w\in W$, called a multiplicity parameter. We write $k_{i}$ for $k_{\alpha_{i}}$ if
$\{\alpha_{1},\cdots,\alpha_{n}\}$ is a basis of simple roots for $R$. It is also clear that there are at most two
$W$-orbits if $R$ is reduced and irreducible. So for convenience we also write $k$ for $k_{1}$ and
$k'$ for $k_{n}$ if $\alpha_{n}$ is not in the $W$-orbit of $\alpha_{1}$. In our situation the root system $R$
of type $A_{n}$ is somehow peculiar, it has only one single $W$-orbit, but we also let $k'$, given in Lemma $\ref{bmap}$,
enter into the parameter $\kappa$, since for type $A_{n}$ there exists a nontrivial $W$-equivariant symmetric bilinear map.
Let
\[
a^{\kappa}: \mathfrak{h}\times\mathfrak{h}\rightarrow\mathbb{C},\quad
b^{\kappa}:\mathfrak{h}\times\mathfrak{h}\rightarrow\mathfrak{h}.
\]
be a $W$-invariant symmetric bilinear form and a $W$-equivariant symmetric bilinear map on $\mathfrak{h}$,
depending on $\kappa$, respectively.

Taking cue from the special hypergeometric functions constructed by Heckman and Opdam \cite{Heckman-Opdam-1,Heckman-Opdam-2}
\cite{Opdam-1,Opdam-2},
we consider for $u,v\in\mathfrak{h}$,
such a second order differential operator on $\mathcal{O}_{H^{\circ}}$ defined by
\[ D^{\kappa}_{u,v}:=\partial_{u}\partial_{v}+\frac{1}{2}\sum\limits_{\alpha>0}
k_{\alpha}\alpha(u)\alpha(v)X_{\alpha}+\partial_{b^{\kappa}(u,v)}+a^{\kappa}(u,v), \]
where the vector fields $X_{\alpha}$'s are defined as
\[   X_{\alpha}:=\frac{e^{\alpha}+1}{e^{\alpha}-1}\partial_{\alpha^{\vee}}.  \]
Notice that $X_{\alpha}$ is invariant under inversion: $X_{-\alpha}=X_{\alpha}$.

It adds to the main linear second order term
a lower order perturbation, consisting of a $W$-equivariant first order term and a $W$-invariant constant.
Notice that $wD^{\kappa}_{u,v}w^{-1}=D^{\kappa}_{wu,wv}$ where $w\in W$.

Inspired by these data, we define connections
$\nabla^{\kappa}=\nabla^{0}+\Omega^{\kappa}$
on the cotangent bundle of $H^{\circ}$ with
$\Omega^{\kappa}\in \mathrm{Hom}(\Omega_{H^{\circ}},\Omega_{H^{\circ}}\otimes\Omega_{H^{\circ}})$
given by
\begin{equation}\label{eq1}
\Omega^{\kappa}: \zeta\in\Omega_{H^{\circ}}\mapsto\frac{1}{2}\sum\limits_{\alpha>0}
k_{\alpha}\zeta(X_{\alpha})d\alpha\otimes d\alpha + (B^{\kappa})^{*}(\zeta).
\end{equation}

Then following the construction in Proposition 2.2 of \cite{Shen}, we define connections
$\tilde{\nabla}^{\kappa}=\tilde{\nabla}^{0}+\tilde{\Omega}^{\kappa}$
on the cotangent bundle of $H^{\circ}\times\mathbb{C}^{\times}$ with
$\tilde{\Omega}^{\kappa}\in \mathrm{Hom}(\Omega_{H^{\circ}\times\mathbb{C}^{\times}},
\Omega_{H^{\circ}\times\mathbb{C}^{\times}}\otimes\Omega_{H^{\circ}\times\mathbb{C}^{\times}})$
given by
\begin{equation}\label{eq2}
\tilde{\Omega}^{\kappa}:
\left\{
\begin{aligned}
\zeta\in\Omega_{H^{\circ}}&\mapsto\frac{1}{2}\sum\limits_{\alpha>0}
k_{\alpha}\zeta(X_{\alpha})d\alpha\otimes d\alpha+(B^{\kappa})^ {*}(\zeta)-\zeta\otimes\frac{dt}{t}-\frac{dt}{t}\otimes\zeta, \\
\frac{dt}{t}\in\Omega_{\mathbb{C}^{\times}} &\mapsto A^{\kappa}-\frac{dt}{t}\otimes\frac{dt}{t}.
\end{aligned}
\right.
\end{equation}
Here $t$ is the coordinate for $\mathbb{C}^{\times}$, and $A^{\kappa}$ resp.
$B^{\kappa}$ denotes the translation invariant
tensor field on $H$ (or $H\times\mathbb{C}^{\times}$) defined by
$a^{\kappa}$ resp. $b^{\kappa}$.

According to ($\ref{eq2}$), we can write
$\tilde{\Omega}^{\kappa}$ explicitly:
\begin{align*}
\tilde{\Omega}^{\kappa}:=&\frac{1}{2}\sum_{\alpha>0}k_{\alpha}d\alpha\otimes d\alpha\otimes X_{\alpha}+(B^{\kappa})^{*}
+c^{\kappa}\sum_{\alpha>0}d\alpha\otimes d\alpha\otimes t\frac{\partial}{\partial t}\\
&-\sum_{\alpha_{i}\in \mathfrak{B}}d\alpha_{i}\otimes\frac{dt}{t}\otimes\partial_{p_{i}}
-\frac{dt}{t}\otimes\frac{dt}{t}\otimes t\frac{\partial}{\partial t}
-\sum_{\alpha_{i}\in \mathfrak{B}}\frac{dt}{t}\otimes d\alpha_{i}\otimes\partial_{p_{i}}.
\end{align*}
Here $c^{\kappa}$ is a constant for each $\kappa$ such that $A^{\kappa}=c^{\kappa}\sum_{\alpha>0}d\alpha\otimes d\alpha$,
$\mathfrak{B}$ is a fundamental system for $R$, and ${p_{i}}$ is the dual basis of $\mathfrak{h}$ to $\alpha_{i}$
such that $\alpha_{i}(p_{j})=\delta^{i}_{j}$ where $\delta^{i}_{j}$ is the Kronecker delta.

Since $\nabla^{\kappa}$ is torsion free: for taking the values in the symmetric tensors,
it is clear that $\tilde{\nabla}^{\kappa}$ is also torsion free.
In \cite{Shen} we prove that
\begin{theorem}\label{thm:flatness}
There exists a bilinear form $a^{\kappa}$ for each $\kappa$ such that $\tilde{\nabla}^{\kappa}$ is flat.
\end{theorem}

\begin{proof}
See Section 2 of \cite{Shen}, where one can also find an explicit form of $a^{\kappa}$ for a given $\kappa$
as follows:
\begin{align*}
&A_{n}: a^{\kappa}(u,v)=\frac{(n+1)}{4}(k^{2}-k'^{2})(u,v);\\
&B_{n}: a^{\kappa}(u,v)=((n-2)k^{2}+kk')(u,v);\\
&C_{n}: a^{\kappa}(u,v)=((n-2)k^{2}+2kk')(u,v);\\
&D_{n}: a^{\kappa}(u,v)=(n-2)k^{2}(u,v);\\
&E_{n}: a^{\kappa}(u,v)=ck^{2}(u,v),\quad c=6,12,30 \; \text{for} \; n=6,7,8;\\
&F_{4}: a^{\kappa}(u,v)=(k+k')(2k+k')(u,v);\\
&G_{2}: a^{\kappa}(u,v)=\frac{3}{4}(k+3k')(k+k')(u,v),
\end{align*}
for which we use the construction of root systems in Bourbaki and take the inner product
$(\cdot,\cdot)$ such
that $(\varepsilon_{i},\varepsilon_{j})=\delta^{i}_{j}$.
\end{proof}

\begin{remark}
Fixing this $a^{\kappa}$, by \cite{Shen} we also know that for every sublattice $L$ of the root lattice $Q$
spanned by elements of $R$,
the `linearized connection' on $\mathfrak{h}-\cup_{\alpha\in R\cap L} \mathfrak{h}_{\alpha}$ defined by the following
$\mathrm{End}(\mathfrak{h})$-valued differential
\[
\Omega_{L}:=\sum_{\alpha\in R\cap L} k_{\alpha}\frac{d\alpha}{\alpha}\otimes\pi_{\alpha}
\]
is flat, where $\pi_{\alpha}\in\mathrm{End}(\mathfrak{h})$ is twice of
the orthogonal projection to $\alpha^{\vee}$ with kernel $\mathfrak{h}_{\alpha}$.
In the meantime each $u_{\alpha}:=k_{\alpha}(\alpha^{\vee}\otimes\alpha)$ is self-adjoint relative to $a^{\kappa}$.
Then $(H,R,\kappa)$ defines a toric analogue of the Dunkl system in the sense of Couwenberg-Heckman-Looijenga
\cite{Couwenberg-Heckman-Looijenga}. Hence the connection $\tilde{\nabla}^{\kappa}$ defined in ($\ref{eq2}$) is usually called
a \emph{toric Dunkl connection}.
\end{remark}

Fixing the bilinear form $a^{\kappa}$ for which $\tilde{\nabla}^{\kappa}$ is flat, we then consider the dual
connections defined on the tangent bundle
of $H^{\circ}\times\mathbb{C}^{\times}$ instead of the connections $\tilde{\nabla}^{\kappa}$  defined on the cotangent bundle of
$H^{\circ}\times\mathbb{C}^{\times}$:
\begin{align*}
-(\tilde{\Omega}^{\kappa})^{*}=&\frac{1}{2}\sum_{\alpha >0}
k_{\alpha}\frac{1+e^{\alpha}}{1-e^{\alpha}}
d\alpha\otimes \partial_{\alpha^{\vee}}\otimes d\alpha
-((B^{\kappa})^{*})'-c^{\kappa}\sum_{\alpha >0}d\alpha\otimes t\frac{\partial}{\partial t}\otimes d\alpha\\
&+\sum_{\alpha_{i}\in \mathfrak{B}}d\alpha_{i}\otimes\partial_{p_{i}}\otimes\frac{dt}{t}
+\sum_{\alpha_{i}\in \mathfrak{B}}\frac{dt}{t}\otimes \partial_{p_{i}}\otimes d\alpha_{i}
+\frac{dt}{t}\otimes t\frac{\partial}{\partial t}\otimes \frac{dt}{t}.
\end{align*}
That the connection form is $-(\tilde{\Omega}^{\kappa})^{*}$ is because the dual connection is characterized by
the property that the pairing between differentials and vector fields is flat.
But in what follows we will still write the dual connection as $\tilde{\nabla}^{\kappa}$
if no confusion would arise.

Since $T_{(p,t)}(H^{\circ}\times\mathbb{C}^{\times})=T_{p}H^{\circ}\oplus T_{t}\mathbb{C}^{\times}$,
we can write a vector field $\tilde{X}$ on $H^{\circ}\times\mathbb{C}^{\times}$ in the following form:
$$\tilde{X}=X(p,t)+\lambda_{1}(p,t)t\frac{\partial}{\partial t},$$
for which $X(p,t)$ is a vector field on $H^{\circ}$ and $\lambda_{1}(p,t)$ is a holomorphic
function depending on both $p$ and $t$. Here we write $\tilde{X}=X+\lambda_{1}t\frac{\partial}{\partial t}$
just for convenience.

Likewise, write $\tilde{Y}=Y+\lambda_{2}t\frac{\partial}{\partial t}$. Inspired by the flat connection $\tilde{\nabla}^{\kappa}$,
we define a product for each $\kappa$ on the
tangent bundle of $H^{\circ}\times\mathbb{C}^{\times}$ by
\begin{align} \label{eq3}
  \tilde{X}\cdot_{\kappa}\tilde{Y}
  :=&\frac{1}{2}\sum_{\alpha >0}
  k_{\alpha}\frac{1+e^{\alpha}}{1-e^{\alpha}}
  \alpha(X) \alpha(Y) \alpha^{\vee}
  -b^{\kappa}(X,Y)-c^{\kappa}\sum_{\alpha >0}\alpha(X) \alpha(Y)t\frac{\partial}{\partial t} \notag \\
&+\sum_{\alpha_{i}\in \mathfrak{B}}\alpha_{i}(X)\lambda_{2}p_{i}
  +\sum_{\alpha_{i}\in \mathfrak{B}}\lambda_{1} p_{i}\alpha_{i}(Y)
  +\lambda_{1}\lambda_{2} t\frac{\partial}{\partial t} \notag \\
=&\frac{1}{2}\sum_{\alpha >0}k_{\alpha}\frac{1+e^{\alpha}}{1-e^{\alpha}}
  \alpha(X) \alpha(Y) \alpha^{\vee}
  -b^{\kappa}(X,Y)-a^{\kappa}(X,Y)t\frac{\partial}{\partial t} \notag \\
&+\lambda_{2}X+\lambda_{1}Y+\lambda_{1}\lambda_{2}t\frac{\partial}{\partial t}.
\end{align}

We already know that $a^{\kappa}$ is a symmetric bilinear form on $\mathfrak{h}$:
\[ a^{\kappa}: \mathfrak{h}\times\mathfrak{h}\rightarrow\mathbb{C}.\]

We can extend $a^{\kappa}$ to be a symmetric bilinear form on $\mathfrak{h}\oplus\mathbb{C}$,
the tangent space of $H^{\circ}\times\mathbb{C}^{\times}$ at $(p,t)$,
by defining
\begin{equation*}
\left\{
\begin{aligned}
&a^{\kappa}(X, t\frac{\partial}{\partial t})=0 \\
&a^{\kappa}(t\frac{\partial}{\partial t}, t\frac{\partial}{\partial t})=-1.
\end{aligned}
\right.
\end{equation*}

Now is $a^{\kappa}(\alpha^{\vee}, \cdot)$ a linear form whose zero set is the hyperplane which
is perpendicular to $\alpha$ and therefore it is proportional to $\alpha$. By
evaluating both sides on $\alpha^{\vee}$ we see that
$$a^{\kappa}(\alpha^{\vee}, \cdot)=\frac{a^{\kappa}(\alpha^{\vee},\alpha^{\vee})}{\alpha(\alpha^{\vee})}\alpha.$$

\begin{remark} \label{rem:identity}
We also notice that
\[(t\frac{\partial}{\partial t})\cdot_{\kappa}\tilde{Y}=Y+\lambda_{2}t\frac{\partial}{\partial t}=\tilde{Y},\]
from which we can see that $t\frac{\partial}{\partial t}$ plays a role of identity in this algebra.
\end{remark}

\begin{theorem} \label{frobenius thm}
The product structure $\cdot_{\kappa}$ defined on $T(H^{\circ}\times\mathbb{C}^{\times})$ by ($\ref{eq3}$)
endows each fiber of $T(H^{\circ}\times\mathbb{C}^{\times})$ with a Frobenius algebra structure.
\end{theorem}

\begin{proof}
In order to see this product structure indeed defines a Frobenius algebra on each fiber of the tangent bundle
of $H^{\circ}\times\mathbb{C}^{\times}$,
we need to verify 3 properties:
\begin{enumerate}
\item [1.] the product is commutative,

\item [2.] the product satisfies the associativity law with respect to the symmetric bilinear form $a^{\kappa}$,
       with this property the trace map can be determined by
       Lemma $\ref{lem:associativity-law}$,

\item [3.] the product is associative.
\end{enumerate}

\noindent $\textbf{1. commutativity of the product}$.

This is quite obvious since the expression for $\tilde{X}\cdot_{\kappa}\tilde{Y}$ is symmetric in
$\{\tilde{X}, \tilde{Y}\}$.

\noindent $\textbf{2. Frobenius condition}$.

Write $\tilde{Z}=Z+\lambda_{3}t\frac{\partial}{\partial t}$, then we have

\begin{align*}
&a^{\kappa}(\tilde{X}\cdot_{\kappa}\tilde{Y},\tilde{Z})\\
=&\frac{1}{2}\sum_{\alpha >0}
k_{\alpha}\frac{1+e^{\alpha}}{1-e^{\alpha}}
      \alpha(X) \alpha(Y) a^{\kappa}(\alpha^{\vee},\tilde{Z})-a^{\kappa}(b^{\kappa}(X,Y),Z)
     -a^{\kappa}(X,Y)a^{\kappa}(t\frac{\partial}{\partial t},\tilde{Z})\\
    &+\lambda_{2}a^{\kappa}(X,\tilde{Z})+\lambda_{1}a^{\kappa}(Y,\tilde{Z})+\lambda_{1}\lambda_{2}a^{\kappa}(t\frac{\partial}{\partial t},\tilde{Z})\\
=&\frac{1}{2}\sum\limits_{\alpha >0}k_{\alpha}\frac{1+e^{\alpha}}{1-e^{\alpha}}
       \cdot\frac{a^{\kappa}(\alpha^{\vee},\alpha^{\vee})}{\alpha(\alpha^{\vee})}\alpha(X) \alpha(Y) \alpha(Z)
       +a^{\kappa}(b^{\kappa}(X,Y),Z)\\
& +\lambda_{3}a^{\kappa}(X,Y)+\lambda_{2}a^{\kappa}(X,Z)+\lambda_{1}a^{\kappa}(Y,Z)-\lambda_{1}\lambda_{2}\lambda_{3}.
\end{align*}
From this, we can see that
\[ a^{\kappa}(\tilde{X}\cdot_{\kappa}\tilde{Y},\tilde{Z})=a^{\kappa}(\tilde{X},\tilde{Y}\cdot_{\kappa}\tilde{Z}),\]
since this expression is fully symmetric in $\lbrace\tilde{X},\tilde{Y},\tilde{Z}\rbrace$. In fact, the symmetry of
$a^{\kappa}(b^{\kappa}(X,Y),Z)$ can be seen from Remark $\ref{rem:cubic-symmetry}$.

\noindent \textbf{3. associativity of the product}.

Let us look at the connection $\tilde{\nabla}^{\kappa}(\mu)$ defined by
\[
\tilde{\nabla}^{\kappa}(\mu)_{\tilde{X}}\tilde{Y}:=\tilde{\nabla}^{0}_{\tilde{X}}\tilde{Y}+\mu\tilde{X}\cdot_{\kappa}
\tilde{Y}.
\]
Written out,
\begin{multline*}
\tilde{\nabla}^{\kappa}(\mu)_{\tilde{X}}\tilde{Y}=
 \tilde{\nabla}^{0}_{\tilde{X}}\tilde{Y}+
 \frac{1}{2}\mu\sum_{\alpha >0}k_{\alpha}\frac{1+e^{\alpha}}{1-e^{\alpha}}
  \alpha(X) \alpha(Y) \alpha^{\vee}
  -\frac{1}{2}\mu k'\sum_{\alpha >0}\alpha(X)\alpha(Y)\alpha'\\
  -\mu c^{\kappa}\sum_{\alpha >0}\alpha(X)\alpha(Y)t\frac{\partial}{\partial t}
  +\mu\lambda_{2}X+\mu\lambda_{1}Y+\mu\lambda_{1}\lambda_{2}t\frac{\partial}{\partial t}.
\end{multline*}
Note that the term $\frac{1}{2}\mu k'\sum_{\alpha >0}\alpha(X)\alpha(Y)\alpha'$ only exists for $A_{n}$ case.

The connection form of $\tilde{\nabla}^{\kappa}(\mu)$ is a holomorphic differential $1$-form on
$H^{\circ}\times\mathbb{C}^{\times}$ taking values in $\mathrm{End}(\mathfrak{h}\oplus\mathbb{C})$. Upon
replacing these endomorphisms, denoted by $\rho_{\alpha}$ or $\rho_{t}$, by their $\mu$ multiplication
$\mu\rho_{\alpha}$ or $\mu\rho_{t}$,
we see that it suffices to prove the flatness of $\tilde{\nabla}^{\kappa}(1)$.
But $\tilde{\nabla}^{\kappa}(1)$ is just $\tilde{\nabla}^{\kappa}$ and we already know that $\tilde{\nabla}^{\kappa}$
is flat by Theorem $\ref{thm:flatness}$, so we can see that $\tilde{\nabla}^{\kappa}(\mu)$ is also flat
for all $\mu\in\mathbb{C}$.
Therefore, the associativity of the product
follows by Proposition $\ref{Frobenius-conditions}$.
\end{proof}

\begin{remark}
In fact, our Frobenius algebra given above includes the Frobenius algebra constructed by Bryan and Gholampour
in \cite{Bryan-Gholampour}
as a special case, which requires $k'=0$ for type $A_{n}$ and $k=k'$ for type $BCFG$. They provided a proof for the
associativity of the product from a point of view of Gromov-Witten theory.
\end{remark}

\begin{corollary}
The Weyl group acts on the tangent bundle by automorphisms. Namely, if we define
$$w(e^{\alpha})=e^{w(\alpha)}$$
for $w\in W$, then for $\tilde{X},\tilde{Y}\in \Gamma(T(H^{\circ}\times\mathbb{C}^{\times}))$,
we have
$$w(\tilde{X}\cdot_{\kappa}\tilde{Y})=w(\tilde{X})\cdot_{\kappa}w(\tilde{Y}).$$
\end{corollary}

\begin{proof}
Let $s_{\beta}$ be the reflection about the hyperplane orthogonal to $\beta$.
By \cite{Bourbaki}, $s_{\beta}$ permute the positive roots other than $\beta$. And since
the terms
$$\frac{1+e^{\alpha}}{1-e^{\alpha}}\partial_{\alpha^{\vee}}\quad \text{and} \quad
\alpha(X)\alpha(Y)$$
remain unchanged under $\alpha\rightarrow -\alpha$, the effect of $s_{\beta}$ to the
formula for $\tilde{X}\cdot_{\kappa}\tilde{Y}$ is to permute the order of the sum:
\begin{align*}
&s_{\beta}(\tilde{X}\cdot_{\kappa}\tilde{Y})\\
=&\frac{1}{2}\sum\limits_{\alpha >0}
        k_{s_{\beta}(\alpha)}\frac{1+e^{s_{\beta}(\alpha)}}{1-e^{s_{\beta}(\alpha)}}s_{\beta}(\alpha)(s_{\beta}X)
        s_{\beta}(\alpha)(s_{\beta}Y)\partial_{s_{\beta}(\alpha^{\vee})}
        -b^{\kappa}(s_{\beta}X,s_{\beta}Y)\\
      &-a^{\kappa}(s_{\beta}X,s_{\beta}Y)s_{\beta}(t\frac{\partial}{\partial t})+s_{\beta}(\lambda_{2}X)
           +s_{\beta}(\lambda_{1}Y)+s_{\beta}(\lambda_{1}\lambda_{2}t\frac{\partial}{\partial t})\\
=&\frac{1}{2}\sum\limits_{\alpha >0}
        k_{\alpha}\frac{1+e^{\alpha}}{1-e^{\alpha}}\alpha(s_{\beta}X)
        \alpha(s_{\beta}Y)\partial_{\alpha^{\vee}}-b^{\kappa}(s_{\beta}X,s_{\beta}Y)\\
      &-a^{\kappa}(s_{\beta}X,s_{\beta}Y)t\frac{\partial}{\partial t}+\lambda_{2}s_{\beta}X
           +\lambda_{1}s_{\beta}Y+\lambda_{1}\lambda_{2}t\frac{\partial}{\partial t}\\
=&(s_{\beta}X+\lambda_{1}t\frac{\partial}{\partial t})\cdot_{\kappa}(s_{\beta}Y+\lambda_{2}t\frac{\partial}{\partial t})\\
=&(s_{\beta}X+s_{\beta}(\lambda_{1}t\frac{\partial}{\partial t}))\cdot_{\kappa}(s_{\beta}Y+s_{\beta}(\lambda_{2}t\frac{\partial}{\partial t}))\\
=&s_{\beta}(\tilde{X})\cdot_{\kappa}s_{\beta}(\tilde{Y})
\end{align*}
since $s_{\beta}(\lambda_{i}t\frac{\partial}{\partial t})=\lambda_{i}t\frac{\partial}{\partial t}$.
Then the corollary follows.
\end{proof}

We thus construct a $W$-invariant fiberwise Frobenius algebra on $H^{\circ}\times\mathbb{C}^{\times}$.
We then have the following theorem.

\begin{theorem}\label{thm:Frobenius-manifolds}
The manifold $H^{\circ}\times\mathbb{C}^{\times}$ endowed with the structure $(\cdot_{\kappa},a^{\kappa},
t\frac{\partial}{\partial t})$ is a Frobenius manifold.
\end{theorem}

\begin{proof}
By Remark $\ref{rem:identity}$, the vector field $t\frac{\partial}{\partial t}$ is the identity of this algebra.
Then we know that $(\cdot_{\kappa},a^{\kappa},
t\frac{\partial}{\partial t})$ endows with a Frobenius algebra on $H^{\circ}\times\mathbb{C}^{\times}$ fiberwisely,
since the trace map $F$ can be determined by the bilinear form $a^{\kappa}$.

We then check the conditions of Definition $\ref{def:Frobenius-manifold}$.
Condition (1) is satisfied since
it is already proved that $\tilde{\nabla}^{\kappa}(\mu)$ is flat for any $\mu\in\mathbb{C}$.
Condition (2) is also clear because the vector field $t\frac{\partial}{\partial t}$ is flat with respect to
the Levi-Civita connection $\tilde{\nabla}^{0}$ of $a^{\kappa}$.
Therefore, $(H^{\circ}\times\mathbb{C}^{\times},\cdot_{\kappa},a^{\kappa},
t\frac{\partial}{\partial t})$ is a Frobenius manifold.
\end{proof}

We also have the dilatation field on $H^{\circ}\times\mathbb{C}^{\times}$ as follows.
\begin{corollary}
Suppose an affine structure on $H^{\circ}\times\mathbb{C}^{\times}$ is given by the torsion free flat connection
$\tilde{\nabla}^{\kappa}$
defined by ($\ref{eq2}$), then the vector field $t\frac{\partial}{\partial t}$ is in fact a
dilatation field on $H^{\circ}\times\mathbb{C}^{\times}$ with factor $\nu=1$.
\end{corollary}

\begin{proof}
It is a straightforward computation. Suppose a local vector field $\tilde{X}$ on $H^{\circ}\times\mathbb{C}^{\times}$ is of the form
$\tilde{X}:=X+\lambda t\frac{\partial}{\partial t}$ where $X$ is a vector field on $H^{\circ}$, we have
\begin{align*}
\tilde{\nabla}_{\tilde{X}}^{\kappa}(t\frac{\partial}{\partial t})&=\tilde{\nabla}_{X+\lambda t\frac{\partial}{\partial t}}^{0}(t\frac{\partial}{\partial t})
-\tilde{\Omega}_{X+\lambda t\frac{\partial}{\partial t}}^{\kappa,*}(t\frac{\partial}{\partial t})\\
&=0-0-0-0+\sum_{\alpha_{i}\in\mathfrak{B}}\alpha_{i}(X)\partial_{p_{i}}+\lambda t\frac{\partial}{\partial t}\\
&=\tilde{X}
\end{align*}
since $t\frac{\partial}{\partial t}$ is flat with respect to $\tilde{\nabla}^{0}$.
\end{proof}

Now let us try to find the (local) potential function for this Frobenius structure. In order to find this
potential function $\Phi$, we require that for $\tilde{X}$, $\tilde{Y}$ and $\tilde{Z}$ being flat vector fields
on $H^{\circ}\times\mathbb{C}^{\times}$, we should have
\begin{align*}
\tilde{\nabla}_{\tilde{X}}^{\kappa}\tilde{\nabla}_{\tilde{Y}}^{\kappa}\tilde{\nabla}_{\tilde{Z}}^{\kappa}\Phi
=&a^{\kappa}(\tilde{X}\cdot\tilde{Y},\tilde{Z})\\
=&\frac{1}{2}\sum_{\alpha>0}k_{\alpha}\frac{1+e^{\alpha}}{1-e^{\alpha}}\cdot
\frac{a^{\kappa}(\alpha^{\vee},\alpha^{\vee})}{\alpha(\alpha^{\vee})}\alpha(X)\alpha(Y)\alpha(Z)
-a^{\kappa}(b^{\kappa}(X,Y),Z)\\
&+\lambda_{1}a^{\kappa}(Y,Z)+\lambda_{2}a^{\kappa}(X,Z)+\lambda_{3}a^{\kappa}(X,Y)
-\lambda_{1}\lambda_{2}\lambda_{3}.
\end{align*}
So let us analyze these terms one by one. For terms $-\lambda_{1}\lambda_{2}\lambda_{3}$ and
$\lambda_{1}a^{\kappa}(Y,Z)+\lambda_{2}a^{\kappa}(X,Z)+\lambda_{3}a^{\kappa}(X,Y)$,
we can easily find their potential functions as follows: $$-\frac{t^{3}}{3!} \quad \text{and} \quad
\frac{t}{2}c^{\kappa}\sum_{\alpha>0}
\alpha^{2},$$ i.e.,
\begin{align*}
&\tilde{\nabla}_{\tilde{X}}^{\kappa}\tilde{\nabla}_{\tilde{Y}}^{\kappa}\tilde{\nabla}_{\tilde{Z}}^{\kappa}(-\frac{t^{3}}{3!})
=-\lambda_{1}\lambda_{2}\lambda_{3}\\
&\tilde{\nabla}_{\tilde{X}}^{\kappa}\tilde{\nabla}_{\tilde{Y}}^{\kappa}\tilde{\nabla}_{\tilde{Z}}^{\kappa}
(\frac{t}{2}c^{\kappa}\sum_{\alpha>0}\alpha^{2})
=\lambda_{1}a^{\kappa}(Y,Z)+\lambda_{2}a^{\kappa}(X,Z)+\lambda_{3}a^{\kappa}(X,Y).
\end{align*}
By the discussion on toric Lauricella case in Section $\ref{sec:toric-Lauricella-manifolds}$,
we can write the term $a^{\kappa}(b^{\kappa}(X,Y),Z)=d^{\kappa}\sum_{\alpha>0}\alpha(X)\alpha(Y)\alpha(Z)$
where $d^{\kappa}$ is a constant when $\kappa$ is given.
So for term $-a^{\kappa}(b^{\kappa}(X,Y),Z)$, we have its potential function:
$$-d^{\kappa}\sum_{\alpha>0}\frac{\alpha^{3}}{3!},$$ i.e.,
\begin{align*}
\tilde{\nabla}_{\tilde{X}}^{\kappa}\tilde{\nabla}_{\tilde{Y}}^{\kappa}\tilde{\nabla}_{\tilde{Z}}^{\kappa}
(-d^{\kappa}\sum_{\alpha>0}\frac{\alpha^{3}}{3!})
=-a^{\kappa}(b^{\kappa}(X,Y),Z).
\end{align*}
Now we have only one term left: $\frac{1}{2}\sum_{\alpha>0}k_{\alpha}\frac{1+e^{\alpha}}{1-e^{\alpha}}\cdot
\frac{a^{\kappa}(\alpha^{\vee},\alpha^{\vee})}{\alpha(\alpha^{\vee})}\alpha(X)\alpha(Y)\alpha(Z)$. It is
not easy to find an explicit potential function for $\frac{1}{2}\cdot\frac{1+e^{\alpha}}{1-e^{\alpha}}$,
but we can always do the Taylor expansion for $\frac{1}{2}\cdot\frac{1+e^{\alpha}}{1-e^{\alpha}}$ and find its
potential series. Let us assume $q(\alpha)$ is a series satisfying
\[
q'''(\alpha)=\frac{1}{2}\cdot\frac{1+e^{\alpha}}{1-e^{\alpha}},
\]
then by the above discussion we have the potential function for this Frobenius structure as follows
\begin{align}\label{eqn:potential-function}
\Phi=-\frac{t^{3}}{3!}+\frac{t}{2}c^{\kappa}\sum_{\alpha>0}\alpha^{2}
+\sum_{\alpha>0}k_{\alpha}\frac{a^{\kappa}(\alpha^{\vee},\alpha^{\vee})}{\alpha(\alpha^{\vee})}q(\alpha)
-d^{\kappa}\sum_{\alpha>0}\frac{\alpha^{3}}{3!},
\end{align}
where $d^{\kappa}=0$ unless for type $A_{n}$. This means we have
\[
\tilde{\nabla}_{\tilde{X}}^{\kappa}\tilde{\nabla}_{\tilde{Y}}^{\kappa}\tilde{\nabla}_{\tilde{Z}}^{\kappa}\Phi
=a^{\kappa}(\tilde{X}\cdot\tilde{Y},\tilde{Z}).
\]

\section{An example: toric Lauricella manifolds}\label{sec:toric-Lauricella-manifolds}

In this section we give an explicit class of Frobenius manifolds, falling into the discussion of
the preceding section. We refer this class of examples as \emph{toric Lauricella
manifolds}. They are called by this name
because their relation to the Lauricella hypergeometric functions \cite{Deligne-Mostow}.

Let $N$ be an index set $\{1,2,\cdots,n+1\}$ and associate to each $i\in N$ a
real number $\mu_{i}\in (0,+\infty)$. Denote the standard basis of $\mathbb{C}^{n+1}$
by $\varepsilon_{1},\cdots,\varepsilon_{n+1}$. We endow $\mathbb{C}^{n+1}$ with a symmetric
bilinear form as $a(z,w):=\sum_{i=1}^{n+1}\mu_{i}z^{i}w^{i}$ for which $z$ is defined by
$z:=\sum z^{i}\varepsilon_{i}$. Let $\mathfrak{h}$ be the quotient of $\mathbb{C}^{n+1}$
by its main diagonal $\Delta_{N}:=\mathbb{C}\sum\varepsilon_{i}$. Since the generator
$\varepsilon_{N}=\sum\varepsilon_{i}$ of the main diagonal has a self-product
$a(\varepsilon_{N},\varepsilon_{N})=\sum \mu_{i}\neq 0$, its
orthogonal complement is nondegenerate. Thus we can often identify $\mathfrak{h}$ with this orthogonal complement,
that is, the hyperplane defined by $\sum\mu_{i}z^{i}=0$,
We take our $\alpha$'s to be the collection $\alpha_{i,j}:=(z_{i}-z_{j})_{i\neq j}$ where $z_{i}$ is the dual
basis of $\varepsilon_{i}$ in $(\mathbb{C}^{n+1})^{*}$. We associate to each $\alpha_{i,j}$
a hyperplane $H_{i,j}$ in $\mathbb{C}^{n+1}$ defined by $\{z^{i}-z^{j}=0\}$, and its orthogonal complement
is spanned by the vector $v_{i,j}:=v_{z_{i}-z_{j}}:=\mu_{j}\varepsilon_{i}-\mu_{i}\varepsilon_{j}$.
It is clear that $v_{i,j}\in\mathfrak{h}$.
We denote by $\mathfrak{h}_{i,j}$ the intersection of $H_{i,j}$ with $\mathfrak{h}$.

We immediately notice that the set $R:=\{\alpha_{i,j}\}$ generates a discrete subgroup
of $\mathfrak{h}^{*}$ whose $\mathbb{R}$-linear span defines a real form $\mathfrak{h}
(\mathbb{R})$ of $\mathfrak{h}$. It's easy to show that $a(v_{i,j},\beta)=0$ for
any $\beta\in \ker(\alpha_{i,j})$. According to \cite{Couwenberg-Heckman-Looijenga},
if $u_{\alpha_{i,j}}$ is the self-adjoint map of
$\mathfrak{h}$ defined by $u_{\alpha_{i,j}}(z)=\alpha_{i,j}(z)v_{i,j}$ (with trace $\mu_{i}+\mu_{j}$),
the system $(\mathfrak{h},\{\mathfrak{h}_{i,j}\},\{\mu_{i}+\mu_{j}\})$ defines a Dunkl system.

As already mentioned in Remark $\ref{rem:k'-for-An}$, there actually exists a nonzero cubic form in this case. Let $\tilde{f}:
\mathbb{C}^{n+1}\rightarrow \mathbb{C}$ be defined by $\tilde{f}(z):=
\sum\mu_{i}(z^{i})^{3}$, and denote by $f:\mathfrak{h}\rightarrow \mathbb{C}$
its restriction to $\mathfrak{h}$. The partial derivative of $\tilde{f}$
with respect to $v_{i,j}$ is $3\mu_{j}\mu_{i}(z_{i}^{2}-z_{j}^{2})$, which
is divisible by $\alpha_{i,j}$.

The symmetric bilinear map $\tilde{b}:\mathbb{C}^{n+1}\times\mathbb{C}^{n+1}
\rightarrow\mathbb{C}^{n+1}$ is defined by $\tilde{b}(\varepsilon_{i},\varepsilon_{j}):=
\delta_{j}^{i}\varepsilon_{i}$. Then the map $b:\mathfrak{h}\times\mathfrak{h}\rightarrow
\mathfrak{h}$ is given as the restriction of $\tilde{b}$ to
$\mathfrak{h}\times\mathfrak{h}$ followed by $\pi:\mathbb{C}^{n+1}\rightarrow\mathfrak{h}$
the orthogonal projection from $\mathbb{C}^{n+1}$ to $\mathfrak{h}$, namely,
$b:=\pi\circ\tilde{b}|_{\mathfrak{h}\times\mathfrak{h}}$. We then
have that $a(\tilde{b}(z,z),z)=\tilde{f}(z)$ and $a(b(z,z),z)=f(z)$. So if we write $\tilde{b}_{i}(z)$ for
$\tilde{b}(\varepsilon_{i},z)$, we can write $\tilde{b}_{i}$ as $\tilde{b}_{i}=
\varepsilon_{i}\otimes z_{i}$. If we write $b_{i,j}(z)$ for $b(v_{i,j},z)$, then
\[ b_{i,j}=\mu_{j}\varepsilon_{i}\otimes z_{i}-\mu_{i}\varepsilon_{j}\otimes z_{j}
-\frac{\mu_{i}\mu_{j}}{\mu_{N}}\varepsilon_{N}\otimes(z_{i}-z_{j}),\]
where $\mu_{N}:=\sum_{i}\mu_{i}$.
If we write $a_{i,j}(z)$ for $a(v_{i,j},z)$, then $a_{i,j}=\mu_{i}\mu_{j}(z_{i}-z_{j})$,
we can verify that $[b_{i,j},b_{k,l}]=-\mu_{N}^{-1}(v_{i,j}\otimes a_{k,l}-v_{k,l}\otimes
a_{i,j})$. Hence we have $[b_{z},b_{w}]=-\mu_{N}^{-1}(z\otimes a_{w}-w\otimes a_{z})$.

\begin{lemma}
The expression $a(b(z_{1},z_{2}),z_{3})$ is symmetric in its arguments if all $\mu_{i}$'s
are equal.
\end{lemma}

\begin{proof}
The lemma is equivalent to saying that for every $z\in\mathfrak{h}$, $b_{z}$ is self-adjoint relative to $a$.

Now we let all $\mu_{i}$ be equal to $1$ in the above example, then the above example becomes the case of a root system
of type $A_{n}$. Since we already know that the dimension of $\mathrm{Hom}(\mathrm{Sym}^{2}\mathfrak{h},\mathfrak{h})^{W}$
is just $1$, then the
$b_{0}=\sum_{\alpha >0}\alpha\otimes\alpha\otimes\alpha'$ given in Lemma $\ref{bmap}$ differs the bilinear map $b$
in the above example just by a scalar. We thus have
$b_{i,j}(z)=z^{i}\varepsilon_{i}-z^{j}\varepsilon_{j}-\frac{1}{n+1}(z^{i}-z^{j})\varepsilon_{N}$.
If $i<j<k$, then
\[ a(b_{i,j}(z),\varepsilon_{j}-\varepsilon_{k})=-z^{j}=a(b_{j,k}(z),\varepsilon_{i}-\varepsilon_{j}); \]
if $i,j,k,l$ are pairwise distinct, then
\[ a(b_{i,j}(z),\varepsilon_{k}-\varepsilon_{l})=0. \]
Since $\{\varepsilon_{i}-\varepsilon_{i+1}\mid i=1,2,\cdots,n\}$ is a basis of $\mathfrak{h}$,
the lemma follows.
\end{proof}

\begin{remark}\label{rem:cubic-symmetry}
Since $a^{\kappa}$ and $b^{\kappa}$ must be a multiple of $a$ and $b$ respectively, the expression
$a^{\kappa}(b^{\kappa}(X,Y),Z)$ is also fully symmetric in its arguments.
\end{remark}

Each $\alpha_{i,j}$ now determines a character $e^{\alpha_{i,j}}$ associated to the exponential map
\[
\exp:\mathfrak{h}\rightarrow H=\mathfrak{h}/2\pi\sqrt{-1} P^{\vee}
\]
to our torus $H$ for which $P^{\vee}$ is the cocharacter lattice relative to the character lattice spanned by $\{\alpha_{i,j}\}$.
Suppose all $\mu_{i}$ being equal now, then the above example becomes our toric case associated to a root system
of type $A_{n}$: $k=\mu_{i}$. Once the symmetric $W$-equivariant bilinear map $b^{\kappa}$ is chosen,
or equivalently, the parameter $k'$ is given.
Then we can define a connection $\tilde{\nabla}^{\kappa}$ on the (co)tangent bundle of $H^{\circ}\times\mathbb{C}^{\times}$
as in ($\ref{eq2}$),
and by Theorem $\ref{thm:flatness}$ there exists a corresponding
bilinear form $a^{\kappa}$,
a multiple of $a$, such that the connection $\tilde{\nabla}^{\kappa}$ is torsion free and flat. Thus by Theorem
$\ref{thm:Frobenius-manifolds}$ we have a class of Frobenius manifolds: toric Lauricella manifolds.
Their fiberwise Frobenius algebra and potential functions are given as in ($\ref{eq3}$) and ($\ref{eqn:potential-function}$)
respectively.

\bibliographystyle{alpha}

\end{document}